\documentclass[review]{elsarticle}

\usepackage{lineno,hyperref}
\modulolinenumbers[5]

\journal{Artificial Intelligence Modeling for Dynamical Problems}
\usepackage{tabularx}
\usepackage{textcomp}
\usepackage{graphics}
\usepackage{epsfig}
\usepackage{amssymb,amsmath}
\usepackage{amsfonts}
\usepackage{color}
\usepackage{epstopdf}
\usepackage{multirow}
\usepackage{physics}
\usepackage{float}
\usepackage{mathtools}
\usepackage{amsthm}
\usepackage{caption}
\usepackage{subcaption}
\newcommand{\vect}[1]{\boldsymbol{#1}}
\newcommand{\params}{\boldsymbol{\theta}}
\newcommand{\reals}{\mathbb{R}}
\newcommand{\adjoin}{\mathbf{a}}
\newcommand{\jacob}[2]{\frac{\partial #1}{\partial #2}}

\newtheorem{theorem}{Theorem}[section]

\newtheorem{definition}[theorem]{Definition}

\catcode`@=11 \@addtoreset{equation}{section}
\renewcommand\theequation{\thesection.\@arabic\c@equation}
\catcode`@=12

\usepackage{tikz}
\usetikzlibrary{positioning, arrows.meta, shapes.geometric, calc, decorations.pathreplacing,shapes,arrows, shadows, fit, backgrounds}

\makeatother




\begin{document}
	
	\begin{frontmatter}
		
		\title{\textbf{Neural Ordinary Differential Equations for Modeling Socio-Economic Dynamics}}

		%\tnotetext[mytitlenote]{Fully documented templates are available in the elsarticle package on \href{http://www.ctan.org/tex-archive/macros/latex/contrib/elsarticle}{CTAN}.}
		
		%% Group authors per affiliation:
		%\author{Elsevier\fnref{myfootnote}}
		%\address{Radarweg 29, Amsterdam}
		%\fntext[myfootnote]{Since 1880.}
	
		%% or include affiliations in footnotes:
		\author[1,3]{Sandeep Kumar Samota\corref{mycorrespondingauthor}}\ead{sandeepksamota@gmail.com}
		\cortext[mycorrespondingauthor]{Corresponding author}

		\author[1,3]{S. Chakraverty}\ead{sne\_chak@yahoo.com}

		\author[2,3]{Narayan Sethi}\ead{nsethinarayan@gmail.com}
		
		\address[1]{Department of Mathematics, National Institute of Technology Rourkela, India.}
		
		\address[2]{Department of Humanities and Social Sciences, National Institute of Technology Rourkela, India.}
		
		\address[3]{Poverty Alleviation Research Centre, National Institute of Technology Rourkela, India.}

		%%%%%%%%%%%%%%%%%%%%%%%%%%%%%%%%%%%%%%%%%%%%%%%%%%%%%%%%%%%%%%%%%%%%%%%%%%%%%%%%%%%%%%%%%%%%%%%%%%%%%%%%%%%%%%%%%%%%%%%%%%

\begin{abstract}
Poverty is a complex dynamic challenge that cannot be adequately captured using predefined differential equations. Nowadays, artificial machine learning (ML) methods have demonstrated significant potential in modelling real-world dynamical systems. Among these, Neural Ordinary Differential Equations (Neural ODEs) have emerged as a powerful, data-driven approach for learning continuous-time dynamics directly from observations. This chapter applies the Neural ODE framework to analyze poverty dynamics in the Indian state of Odisha. Specifically, we utilize time-series data from 2007 to 2020 on the key indicators of economic development and poverty reduction. Within the Neural ODE architecture, the temporal gradient of the system is represented by a multi-layer perceptron (MLP). The obtained neural dynamical system is integrated using a numerical ODE solver to obtain the trajectory of over time. In backpropagation, the adjoint sensitivity method is utilized for gradient computation during training to facilitate effective backpropagation through the ODE solver. The trained Neural ODE model reproduces the observed data with high accuracy. This demonstrates the capability of Neural ODE to capture the dynamics of the poverty indicator of concrete-structured households. The obtained results show that ML methods, such as Neural ODEs, can serve as effective tools for modeling socioeconomic transitions. It can provide policymakers with reliable projections, supporting more informed and effective decision-making for poverty alleviation.
\end{abstract}

%%%%%%%%%%%%%%%%%%%%%%%%%%%%%%%%%%%%%%%%%%%%%%%%%%%%%%%%%%%%%%%%%%%%%%%%%%%%%%%%%%%%%%%%%%%%%%%%%%%%%%%%%%%%%%%%%%%%%%%%%%


\end{frontmatter}

%\linenumbers

%%%%%%%%%%%%%%%%%%%%%%%%%%%%%%%%%%%%%%%%%%%%%%%%%%%%%%%%%%%%%%%%%%%%%%%%%%%%%%%%%%%%%%%%%%%%%%%%%%%%%%%%%%%%%%%%%%%%%

\section{Introduction}

Poverty is a dynamical process that evolves over time under the joint influence of household characteristics, policy interventions, macroeconomic shocks, and structural transformations \cite{carr2009evaluating}. These drivers interact in a nonlinear manner and often exhibit time-varying effects, feedback loops, and delayed responses. Consequently, empirical poverty indicators (e.g., consumption deprivation, multidimensional indices, or asset-based proxies) typically display trajectories that are difficult to represent with a fixed, hand-crafted system of differential equations \cite{kumar2025assessing, kumar2023multidimensional}.

Socioeconomic modeling has traditionally relied on econometric and system-based approaches, including time-series and panel methods, structural models, and mechanistic dynamical systems \cite{deaton1985panel,pesaran2015time,wu2015economic}. Such models are attractive because the assumed structure can be interpreted in terms of economic mechanisms and policy levers; however, they generally require specifying functional forms and state interactions. This becomes challenging when the governing relations are unknown, heterogeneous across regions, or non-stationary due to institutional and technological change. In many poverty settings, data are also sparse, irregularly sampled, or incomplete, which complicates purely discrete-time modeling and motivates methods that can naturally operate in continuous time.

In parallel, machine learning (ML) has been increasingly used for socioeconomic measurement and prediction \cite{kumar2025comparative}. Supervised learning has been applied to poverty mapping and small-area estimation using high-dimensional covariates such as satellite imagery, mobile-phone metadata, and administrative records \cite{jean2016combining,steele2017mapping,khoun2025mapping}. Further, sequence models (e.g., RNNs, GRUs, and LSTMs) have been employed to capture temporal dependence in longitudinal socioeconomic indicators \cite{mcdermott2019bayesian, emshagin2022short}. Furthermore, ML has applications in solving different types of differential equations using physics-informed neural networks (PINNs) and their subsequent variants \cite{sahu2026physics,raissi2019physics,kumar2023physics}. These approaches can learn complex nonlinear relationships directly from data, but many ML pipelines remain fundamentally discrete-time: predictions are tied to fixed observation grids, and handling missingness or irregular sampling often requires ad hoc imputation or resampling. Moreover, for policy analysis, it is often desirable to generate smooth trajectories and counterfactual forecasts at arbitrary time points rather than only at the times present in the training data.

Neural Ordinary Differential Equations (Neural ODEs) provide a principled bridge between continuous-time dynamical systems and deep learning. The core idea, introduced by Chen et al.~\cite{chen2018neural}, is to replace a finite sequence of residual network layers by a continuous-time flow defined through an ODE, $\dot{x}(t)=f_{\theta}(x(t),t)$, where the vector field $f_{\theta}$ is parameterized by a neural network. This formulation makes the ``depth'' of the model continuous and shifts the computation of the forward pass to a numerical ODE solver, enabling evaluation at arbitrary time points. Efficient training is enabled by differentiating through the solver, most commonly via the adjoint sensitivity method, which computes gradients by solving a related backward-in-time ODE.

Subsequent work extended Neural ODEs to latent-variable and generative settings. Latent ODE models, introduced by Rubanova et al.~\cite{rubanova2019latent}, combine a variational encoder (often an RNN/GRU) with a neural ODE defined in latent space, allowing irregularly sampled multivariate time series to be modeled as continuous-time latent trajectories with a probabilistic observation model. Related developments include continuous normalizing flows for density modeling \cite{grathwohl2018ffjord} and augmented/controlled formulations that improve expressiveness and numerical behavior \cite{kidger2020neural,kidger2022neural}. These latent-space Neural ODE formulations are particularly relevant for socioeconomic applications, where the observed indicators can be noisy proxies of underlying latent development states and where measurement times may be uneven across regions or years.

In this chapter, a district-embedded latent Neural ODE architecture is applied to model poverty dynamics in Odisha, India, using a time-series proxy of economic well-being during 2007--2020. District embeddings are used to capture persistent heterogeneity, while a GRU-based encoder summarizes the observed sequence into an initial latent state that is evolved in continuous time by the neural ODE. After reconstructing missing observations, the learned dynamics provide smooth trajectories for interpolation and enable forecasting to future years, supporting evidence-based planning.

The remainder of the chapter is organized as follows. Section~\ref{sec:prelimiries} introduces methodological preliminaries, including the relationship between residual networks and ODE solvers, the Neural ODE formulation, and training via adjoint sensitivity. Section~\ref{sec:neural_odes_for_poverty_dynamics} describe the data, the model architecture and training procedure. Section~\ref{sec:results_and_discussions} describes the empirical results and forecast summaries, and discussions. Finally, Section~\ref{sec:conclusion} concludes the chapter.

\section{Preliminaries}\label{sec:prelimiries}

This section discusses the main concepts required to formulate and train Neural ODE model for socioeconomic time-series.

\subsection{Residual Networks (ResNets)}

To understand the conceptual leap of Neural ODEs, it is important to first understand the ResNet architecture on which they generalize. A standard feedforward neural network (SFNN) transforms a hidden state $\vect{h}_t$ at layer $t$ to a new state $\vect{h}_{t+1}$ at layer $t+1$ using a non-linear activation function $g$ as
\begin{equation}
	\vect{h}_{t+1} = g(\vect{h}_t, \params_t),
\end{equation}
where $\params_t$ are the parameters (weights and biases) of the $t$-th layer. As the number of layers $N$ grows, training becomes extremely difficult due to the vanishing or exploding gradient problem.

The ResNets modifies this update rule by adding a skip connection (or identity mapping) as \cite{he2016deep}
\begin{equation}
	\vect{h}_{t+1} = \vect{h}_t + f(\vect{h}_t, \params_t).
	\label{eq:resnet_update}
\end{equation}
Here, the network $f$ learns the residual (the change) to be added to the input $\vect{h}_t$. This formulation dramatically improves gradient flow during backpropagation, as the gradient can pass directly through the identity connection.

% ----------------------------------------------------------------------
\subsection{The Continuous-Time Analogy}
% ---------------------------------------------------------------------
The ResNet update rule presented in Eq. \eqref{eq:resnet_update} exhibits a substantial formal similarity to a numerical method, viz., the Euler method for solving ODEs. Consider a general initial value problem (IVP) defined by
\begin{equation}
	\frac{d\vect{z}(t)}{dt} = f(\vect{z}(t), t, \params)
	\quad \text{with initial condition} \quad \vect{z}(t_0) = \vect{z}_0.
	\label{eq:ode_def}
\end{equation}
Here, $\vect{z}(t) \in \reals^d$ is the state of a system at time $t$, and $f$ is a function (parameterized by $\params$) that describes the dynamics of $\vect{z}$.

The simplest numerical method to approximate the solution $\vect{z}(t_1)$ is the forward Euler method. It discretizes time into steps of size $\Delta t$ and approximates the derivative as \cite{chien2022learning}
\begin{equation*}
	\frac{\vect{z}(t + \Delta t) - \vect{z}(t)}{\Delta t} \approx f(\vect{z}(t), t, \params).
\end{equation*}
Rearranging this gives the update rule
\begin{equation}
	\vect{z}(t + \Delta t) = \vect{z}(t) + f(\vect{z}(t), t, \params) \Delta t.
	\label{eq:euler_update}
\end{equation}

By comparing Equation (\ref{eq:resnet_update}) and (\ref{eq:euler_update}), the connection is clear as follows:
\begin{itemize}
	\item A ResNet block is equivalent to one step of the forward Euler method with a step size of $\Delta t = 1$.
	\item The layer index $t$ in a ResNet is similar to the time $t$ in the ODE.
	\item The residual function $f(\vect{h}_t, \params_t)$ in a ResNet is similar to the dynamics function $f(\vect{z}(t), t, \params)$ of the ODE.
\end{itemize}

Therefore, the Neural ODE takes this similarity to its logical conclusion that instead of a discrete sequence of layers, one can define a continuous transformation. This continuous-time viewpoint is attractive for poverty dynamics because observations may be irregularly spaced, and policy analysis often requires smooth trajectories and forecasts at arbitrary time points. We can learn the dynamics function $f$ itself, and the layers are created dynamically by the numerical solver.

\subsection{Neural Ordinary Differential Equations}

\begin{definition}[Neural ODE \cite{chen2018neural}]
	A Neural ODE is a continuous-time depth neural network model. Suppose at an initial time $t_0$,
	\begin{equation}
		\vect{z}(t_0) = \vect{z}_0
	\end{equation}
    is the initial state of the hidden dynamics $\vect{z}(t)$ for a given input $\vect{z}_0 \in \reals^d$. The hidden state then evolves according to an ODE, 
	\begin{equation}
		\frac{d\vect{z}(t)}{dt} = f_{\params}(\vect{z}(t), t, \params), \quad t\in [t_0,t_1]
	\end{equation}
    parameterized by a neural network, $f$, with parameters $\params$.
\end{definition}

Then the output of the Neural ODEs system is the state $\vect{z}(t_1)$ at a chosen final time $t_1$. This output is the solution to the IVP and it is found by integrating the dynamics
\begin{equation}\label{eq:forward_integral}
	\vect{z}(t_1) = \vect{z}(t_0) + \int_{t_0}^{t_1} f_\theta(\vect{z}(t), t, \params) dt.
\end{equation}

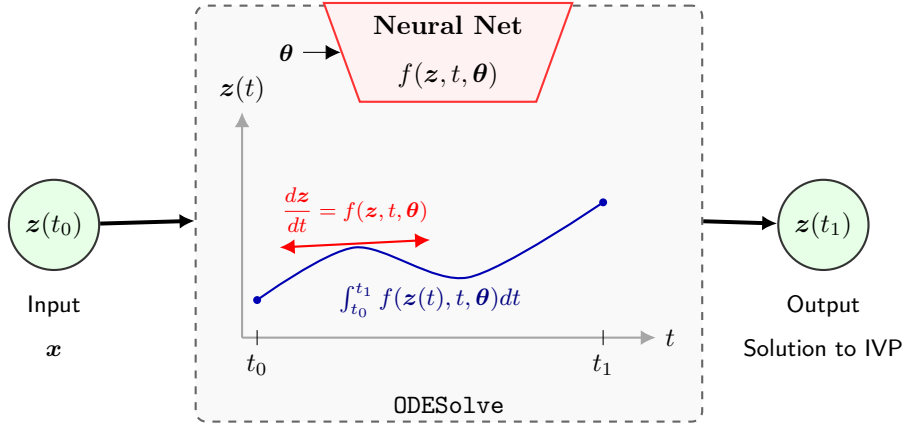
\begin{figure}[ht]
    \centering
    \begin{tikzpicture}[
        node distance=1.5cm and 2cm,
        >={Latex[width=2mm,length=2mm]},
        process/.style={rectangle, draw=black!80, fill=blue!5, thick, minimum width=3cm, minimum height=1.5cm, align=center, rounded corners},
        net/.style={trapezium, trapezium angle=70, draw=red!80, fill=red!5, thick, minimum width=2.5cm, align=center, inner sep=5pt},
        state/.style={circle, draw=black!80, fill=green!10, thick, minimum size=1.2cm, inner sep=0pt},
        label text/.style={font=\small\sffamily, align=center},
        axis/.style={->, thick, draw=gray!70}
    ]

    % --- 1. Initial State ---
    \node[state] (z0) {$\vect{z}(t_0)$};
    \node[below=0.2cm of z0, label text] {Input\\$\vect{x}$};

    % --- 2. The ODE Solver Box (Container) ---
    % We define coordinates to draw the box around later
    \node[right=4.5cm of z0] (center_helper) {};
    
    % --- 3. Internal Visualization of ODE Solver (Trajectory) ---
    % Draw a small plot representing the continuous evolution inside the solver
    \begin{scope}[shift={(2.5, -0.5)}]
        % Axis
        \draw[axis] (0,-1) -- (5.5,-1) node[right] {$t$};
        \draw[axis] (0,-1) -- (0,2.0) node[above] {$\vect{z}(t)$};
        
        % Ticks
        \draw (0.2, -0.9) -- (0.2, -1.1) node[below] {$t_0$};
        \draw (4.8, -0.9) -- (4.8, -1.1) node[below] {$t_1$};
        
        % The Trajectory Curve (The Integral)
        \draw[thick, blue!70!black, smooth] plot coordinates {(0.2,-0.5) (1.5, 0.2) (3.0, -0.2) (4.8, 0.8)};
        \coordinate (t_mid) at (2.2, 1.0); % A point on the curve to tap into
        \fill[blue!70!black] (0.2,-0.5) circle (1.5pt);
        \fill[blue!70!black] (4.8,0.8) circle (1.5pt);
        
        % Tangent Vector (Derivative)
        \draw[<->, red, thick] (0.5, 0.2) -- (2.5, 0.3);
        \node[font=\footnotesize, text=red] at (1.5,0.7) {$\dfrac{d\vect{z}}{dt}=f(\vect{z}, t, \params)$};
        % Label for the integral
        \node[font=\small, text=blue!50!black] at (2.5,-0.5) {$\int_{t_0}^{t_1} f(\vect{z}(t), t, \params) dt$};
    \end{scope}

    % --- 4. The Neural Network f ---
    % Placed above the trajectory to show it parameterizes the derivative
    \node[net, above=1.5cm of center_helper, shape border rotate=180] (nnet) {\textbf{Neural Net} \\ $f(\vect{z}, t, \params)$};
    
    % Connections for the Network
    %\draw[->, dashed, thick, draw=gray] (5.0, 0.5) -- ++(0, 1.5) -| (nnet.east) node[pos=0.2, right, font=\footnotesize] {Current State $\vect{z}(t)$};
    %\draw[->, dashed, thick, draw=red] (nnet.south) -- (5.0, 1.2) node[midway, right, font=\footnotesize] {Dynamics};

    % --- 5. Output State ---
    \node[state, right=9cm of z0] (z1) {$\vect{z}(t_1)$};
    \node[below=0.2cm of z1, label text] {Output\\Solution to IVP};

    % --- 6. Drawing the Solver Box and Main Flow Arrows ---
    
    % Background box for ODESolve
    \begin{scope}[on background layer]
        \node[draw=black!60, fill=gray!5, dashed, thick, rounded corners, fit={(2.0,-2.5) (8.5, 2.8)}, label={[anchor=south]south:\texttt{ODESolve}}] (solver_box) {};
    \end{scope}

    % Main Arrows
    \draw[->, line width=1.5pt] (z0) -- (solver_box);
    \draw[->, line width=1.5pt] (solver_box) -- (z1);

    % --- 7. Annotations ---
    % Parameters
    \node[left=0.5cm of nnet, font=\small] (theta) {$\params$};
    \draw[->] (theta) -- (nnet);
    
    \end{tikzpicture}
    \caption{Visual representation of a Neural ODE. The input $\vect{z}(t_0)$ is evolved continuously to $\vect{z}(t_1)$ by an ODE solver. The dynamics of the evolution (the derivative) are defined by the neural network $f$, parameterized by $\theta$.}
    \label{fig:neural_ode}
\end{figure}

The function $f_{\params}$ can be any neural network architecture that maps the current state $\vect{z}(t)$ and time $t$ to the derivative $d\vect{z}/dt$. The schematic diagram of Neural ODE is shown in Fig. \ref{fig:neural_ode}.

The model therefore learns the governing dynamics from data rather than requiring an explicit, hand-crafted functional form. Depending on the application, $\vect{z}(t)$ may represent the observed quantity directly (e.g., the proportion of concrete-structured households) or a latent state that is mapped to observations through a readout layer.

\subsection{Numerical integration}

Because $f_{\params}$ is nonlinear, closed-form solutions for the integral in Eq. \eqref{eq:forward_integral} are typically unavailable. Therefore, a numerically established black-box ODE solver is utilized to solve this integral. Hence, the integral in Eq. \eqref{eq:forward_integral} is solved as
\begin{equation}
	\vect{z}(t_i) = \texttt{ODESolve}(\vect{z}(t_0), f_{\params}, t_0, t_i).
\end{equation}
This numerically established \texttt{ODESolve} function is the core of this model. Unlike the fixed steps of the Euler method, these types of modern solvers are far more sophisticated. They use higher-order methods like the Dormand-Prince (RKDP) or Adams-Moulton methods, which are significantly more accurate than the forward Euler method. Furthermore, the numerical solver continuously estimates its own local truncation error. If the error is too high (e.g., the dynamics are changing rapidly), it automatically decreases its step size to maintain precision. If the error is low (dynamics are smooth), it increases its step size to be more efficient. 

This adaptive computation means the Neural ODE implicitly adjusts its own depth (number of function evaluations) based on the inputs and the complexity of the transformation. The precision of the model is controlled by a user-defined tolerance, not by a fixed number of layers. In practice, the solver can be evaluated exactly at the available measurement times, which avoids forcing the data into a fixed discrete-time grid.

\subsection{Training objective and gradients}

Assume observations $\{\tilde{\vect{z}}_{i}\}$ are available at times $\{t_{i}\}_{i=1}^{N}$. Parameters are estimated by minimizing a data misfit loss, for example the mean squared error (MSE):
\begin{equation}
\mathcal{L}(\theta)=\frac{1}{N}\sum_{i=1}^{N}\|\vect{z}_{i}-\tilde{\vect{z}}_{i}\|^{2}.
\end{equation}
Training the model requires computing the gradient of a scalar loss function $\mathcal{L}$ with respect to the parameters, $\frac{d\mathcal{L}}{d\params}$. The naive backpropagation approach applys the chain rule through all the discrete operations of the \texttt{ODESolve} function. This is known as backpropagation-through-time (BPTT) if we unroll the steps of ODESolver.

This approach is computationally and memory-prohibitive because:
\begin{enumerate}
	\item \textit{Memory Cost:} It would require storing all intermediate states $\vect{z}(t_i)$ and internal states of the solver computed during the forward pass. This memory cost scales up linearly with the number of steps taken by the solver, which can be in thousands.
	\item \textit{Computational Cost:} The backpropagation path is tied directly to the forward solver, which may not be the most efficient path.
\end{enumerate}
Therefore, this approach completely negates the primary advantage of Neural ODEs.

\subsection{The Adjoint Sensitivity Method}
The adjoint sensitivity method provides an efficient alternative by computing gradients through a backward-in-time sensitivity equation.
It is a technique from optimal control theory \cite{pontryagin2018mathematical}. Instead of backpropagation through the solver, this method allows computing all required gradients by solving a second augmented ODE backward in time. This method provides $\frac{d\mathcal{L}}{d\vect{z}(t)}$ and $\frac{d\mathcal{L}}{d\params}$ with $O(1)$ memory complexity. In the following paragraphs, the adjoint sensitivity method is derived step-by-step \cite{chen2018neural}.

Define the adjoint variable $\adjoin(t)$ as the sensitivity of the loss with respect to the state $\vect{z}(t)$,
\begin{equation}
 \adjoin(t) := \frac{\partial \mathcal{L}}{\partial \vect{z}(t)}\in\mathbb{R}^{d}.
\end{equation}

The $\adjoin(t)$ is treated as a row vector, consistent with its role as a co-vector (a linear map from state perturbations to loss perturbations).

Next goal is to find the dynamics of this adjoint state, i.e., $\frac{d\adjoin(t)}{dt}$. Starting from the initial condition for this backward-in-time problem, which is at time $t_1$ defined as
\begin{equation}
	\adjoin(t_1) = \frac{d\mathcal{L}}{d\vect{z}(t_1)} = \jacob{\mathcal{L}(\vect{z}(t_1))}{\vect{z}(t_1)}.
\end{equation}
This is simply the gradient of the loss function with respect to output of the model, which is the standard starting point for backpropagation.

Next, changes in $\adjoin(t)$ are computed while moving backward in time. Let us consider the relationship between $\vect{z}(t)$ and $\vect{z}(t-\Delta t)$ as
\begin{equation*}
	\vect{z}(t) \approx \vect{z}(t-\Delta t) + f(\vect{z}(t-\Delta t), t-\Delta t, \params) \Delta t.
\end{equation*}
Therefore, using the chain rule, $\adjoin(t)$ can be related to $\adjoin(t-\Delta t)$ as
\begin{align*}
	\adjoin(t-\Delta t) &= \frac{d\mathcal{L}}{d\vect{z}(t-\Delta t)} \\
	&= \frac{d\mathcal{L}}{d\vect{z}(t)} \jacob{\vect{z}(t)}{\vect{z}(t-\Delta t)} \\
	&= \adjoin(t) \jacob{\left( \vect{z}(t-\Delta t) + f(\vect{z}(t-\Delta t), t-\Delta t, \params) \Delta t \right)}{\vect{z}(t-\Delta t)} \\
	&\approx \adjoin(t) \left( \mathbf{I} + \jacob{f(\vect{z}(t-\Delta t), t-\Delta t, \params)}{\vect{z}(t-\Delta t)} \Delta t \right).
\end{align*}
Here, $\mathbf{I}$ is the identity matrix, and $\jacob{f}{\vect{z}}$ is the Jacobian matrix of the dynamics function $f$ with respect to its state input $\vect{z}$.

Rearranging the terms yields
\begin{equation*}
	\adjoin(t-\Delta t) - \adjoin(t) \approx  \adjoin(t) \jacob{f}{\vect{z}}   \Delta t.
\end{equation*}
Dividing by $-\Delta t$, we get
\begin{equation*}
	\frac{\adjoin(t) - \adjoin(t-\Delta t)}{\Delta t} \approx - \adjoin(t) \jacob{f}{\vect{z}}.
\end{equation*}
Now, taking the limit as $\Delta t \to 0$, the ODE that governs the adjoint state is given by
\begin{equation}
	\frac{d\adjoin(t)}{dt} = - \adjoin(t)  \jacob{f_{\params}(\vect{z}(t), t)}{\vect{z}}.
	\label{eq:adjoint_dynamics}
\end{equation}

% For times $t$ where the loss has no direct dependence on $x(t)$ (i.e., between observation points), the adjoint evolves according to the backward ODE
% \begin{equation}
% \frac{d a(t)}{dt} = -\left(\frac{\partial f_{\theta}(x(t),t)}{\partial x}\right)^{\!\top} a(t).
% \label{eq:adjoint_ode}
% \end{equation}

The obtained Eq. \eqref{eq:adjoint_dynamics} is the first component of the backward ODE. It describes how the loss sensitivity flows backward through the continuous dynamics. Notice that its evaluation at time $t$ depends on both $\vect{z}(t)$ and $\adjoin(t)$. The term $\adjoin(t) \jacob{f}{\vect{z}}$ is a vector-Jacobian product (VJP), which is the fundamental operation of reverse-mode automatic differentiation (AD) \cite{paszke2017automatic}.

Next, the gradient $\frac{d\mathcal{L}}{d\params}$ is required. The parameters $\params$ influence the loss $\mathcal{L}$ by affecting the entire trajectory $\vect{z}(t)$ from $t_0$ to $t_1$. The total change in $\mathcal{L}$ due to a change in $\params$ is the integral of the sensitivity contributions over the entire time interval, given by
\begin{equation*}
	\frac{d\mathcal{L}}{d\params} = \int_{t_0}^{t_1} \frac{d\mathcal{L}}{d\vect{z}(t)} \jacob{\vect{z}(t)}{\params} dt.
\end{equation*}
This is not helpful yet. Consider the effect of $\params$ on the loss at a specific time $t$. A change in $\params$ changes the dynamics $f$, which in turn changes the state $\vect{z}(t)$, which in turn changes the final loss $\mathcal{L}$.

Therefore, the total gradient can also be written as an integral over the time-dependent influence as
\begin{equation*}
	\frac{d\mathcal{L}}{d\params} = \int_{t_0}^{t_1} \frac{\partial \mathcal{L}}{\partial f_{\params}(\vect{z}(t), t)} \jacob{f_{\params}(\vect{z}(t), t)}{\params} dt.
\end{equation*}

The obvious question here is how a change in the function value $f$ at time $t$ affects the final loss $\mathcal{L}$. A small perturbation $\delta f$ at time $t$ over a duration $dt$ perturbs the state $\vect{z}(t)$ by $\delta \vect{z}(t) = \delta f \cdot dt$. The effect of this state perturbation on the final loss is $\adjoin(t) \cdot \delta \vect{z}(t)$.

Therefore, the sensitivity of the loss to the function $f$ is exactly the adjoint state
\begin{equation*}
	\frac{\partial \mathcal{L}}{\partial f(\vect{z}(t), t, \params)} = \adjoin(t).
\end{equation*}
Substituting this back into above integral yields
\begin{equation}
	\frac{d\mathcal{L}}{d\params} = \int_{t_0}^{t_1} \adjoin(t) \jacob{f(\vect{z}(t), t, \params)}{\params} dt.
	\label{eq:param_grad_integral}
\end{equation}

This integral can be computed by defining a new state variable, $\vect{g}(t) = \frac{d\mathcal{L}}{d\params}(t)$, whose dynamics are
\begin{equation}
	\frac{d\vect{g}(t)}{dt} = \adjoin(t) \jacob{f(\vect{z}(t), t, \params)}{\params}.
	\label{eq:param_grad_dynamics}
\end{equation}

Integrating this from $t_0$ to $t_1$ yields the total gradient. Alternatively, it can be integrated backward from $t_1$ to $t_0$ by flipping the sign \cite{oh2025comprehensive}
\begin{equation}
	\frac{d\vect{g}(t)}{dt} = -\adjoin(t) \jacob{f(\vect{z}(t), t, \params)}{\params},
	\label{eq:param_grad_dynamics_reverse}
\end{equation}
with the initial condition $\vect{g}(t_1) = \vect{0}$, because the gradient integral starts at zero. Moreover, the final value $\vect{g}(t_0)$ gives the desired gradient $\frac{d\mathcal{L}}{d\params}$.

At this stage, dynamics for $\adjoin(t)$ in Eq. (\ref{eq:adjoint_dynamics}) and $\frac{d\mathcal{L}}{d\params}$ in Eq. (\ref{eq:param_grad_dynamics_reverse}) are available. However, both of these dynamics depend on the original state $\vect{z}(t)$. The entire forward trajectory of $\vect{z}(t)$ could be stored, but this would dethrone the $O(1)$ memory objective. The innovative solution is to recompute the $\vect{z}(t)$ trajectory backward simultaneously with solving for the gradients. This can be done by defining a single, augmented state vector $S(t)$ that concatenates all the states needed, which is given as follows:
\begin{equation}
	S(t) = 
	\begin{bmatrix}
		\vect{z}(t) \\
		\adjoin(t) \\
		\vect{g}(t)
	\end{bmatrix}
	=
	\begin{bmatrix}
		\vect{z}(t) \\
		\frac{d\mathcal{L}}{d\vect{z}(t)} \\
		\frac{d\mathcal{L}}{d\params}(t)
	\end{bmatrix}.
\end{equation}

The dynamics $\frac{dS(t)}{dt}$ of this augmented state can be given by combining the three ODEs derived in Eq. (\ref{eq:ode_def}), (\ref{eq:adjoint_dynamics}), and (\ref{eq:param_grad_dynamics_reverse}) as
\begin{equation}
	\frac{dS(t)}{dt} = \frac{d}{dt} 
	\begin{bmatrix}
		\vect{z}(t) \\
		\adjoin(t) \\
		\vect{g}(t)
	\end{bmatrix}
	=
	\begin{bmatrix}
		f(\vect{z}(t), t, \params) \\
		-\adjoin(t) \jacob{f}{\vect{z}} \\
		-\adjoin(t) \jacob{f}{\params}
	\end{bmatrix}.
	\label{eq:augmented_dynamics}
\end{equation}

To train a Neural ODE, a forward pass is first performed by solving Eq. (\ref{eq:ode_def}) from $t_0$ to $t_1$ to obtain $\vect{z}(t_1)$. All gradients are then computed by solving the augmented ODE defined in Eq. (\ref{eq:augmented_dynamics}) backward in time from $t_1$ to $t_0$.
	
The initial conditions for this backward dynamics to solve at time $t_1$ are given by
\begin{itemize}
	\item $S_z(t_1) = \vect{z}(t_1)$ (the output from the forward pass)
	\item $S_a(t_1) = \adjoin(t_1) = \frac{d\mathcal{L}}{d\vect{z}(t_1)}$ (the gradient from the loss function)
	\item $S_g(t_1) = \vect{g}(t_1) = \vect{0}$ (the gradient integral starts at zero)
\end{itemize}
Again, calling the \texttt{ODESolve} gives
\begin{equation}
	\begin{bmatrix}
		\vect{z}(t_0) \\
		\adjoin(t_0) \\
		\vect{g}(t_0)
	\end{bmatrix}
	=
	\texttt{ODESolve}\left( 
	\begin{bmatrix}
		\vect{z}(t_1) \\
		\adjoin(t_1) \\
		\vect{0}
	\end{bmatrix},
	\frac{dS(t)}{dt}, t_1, t_0, \params
	\right).
\end{equation}
Thus, the final state of the solver at time $t_0$ provides all the required gradients
\begin{itemize}
	\item $\vect{g}(t_0) = \int_{t_1}^{t_0} -\adjoin(t) \jacob{f}{\params} dt = \int_{t_0}^{t_1} \adjoin(t) \jacob{f}{\params} dt = \frac{d\mathcal{L}}{d\params}$
	\item $\adjoin(t_0) = \frac{d\mathcal{L}}{d\vect{z}(t_0)} = \frac{d\mathcal{L}}{d\vect{x}}$ (the gradient with respect to the input $\vect{x}$)
\end{itemize}

This single backward call to \texttt{ODESolve} efficiently computes all parameter gradients and the input gradient while only storing the single augmented state $S(t)$. The memory cost is constant with respect to the integration time or number of solver steps. 

Since all required gradients have been obtained, the parameters $\params$ can now be updated using any optimization technique, such as GD or Adam, to learn the most accurate dynamics $\frac{d \vect{z}(t)}{dt}$.

\section{Neural ODEs for Poverty Dynamics}\label{sec:neural_odes_for_poverty_dynamics}

This section outlines a practical workflow for using Neural ODEs to model poverty dynamics from time-indexed socioeconomic indicators.

\subsection{Data Representation}
The model utilizes a time series dataset having various socio-economic indicators across different geographical districts of Odisha state of India over several time points. In this regard, data is collected for year 2007 from District Level Household and Facility Survey (DLHS)-3 \cite{dlhs3_dataset} for year 2015 and 2020 from National Family Health Survey (NFHS)-4 and NFHS-5 \cite{dhs_dataset}. Input data for the model, denoted as $\mathbf{X}$, is structured as a three-dimensional tensor of shape $(N_d, N_t, N_i)$. Specifically, the dataset comprises: $N_d = 30$ unique districts, $N_t = 3$ observed time points: 2007, 2015, and 2020, and $N_i = 6$ distinct indicators. For model training, these calendar years are normalized to a $[0, 1]$ range using a linear transformation, where 2007 maps to 0 and 2020 maps to 1. This normalization facilitates stable training of the continuous-time ODE model.

All indicator values are treated as proportions, assumed to be within the $[0, 1]$ range. The six indicators included in the dataset are:

\begin{enumerate}
    \item Have Access to toilet facility
    \item Use piped drinking water
    \item Use LPG for cooking
    \item Live in a pucca house
    \item Have Electricity connection
    \item Education Secondary or higher (age 15-49)
\end{enumerate}

This structured data representation allows the Neural ODE model to capture complex, district-specific temporal dynamics for each indicator.

\subsection{District Embeddings}
Each district $d$ is associated with a unique, learnable low-dimensional embedding vector, denoted as $\mathbf{e}_d \in \mathbb{R}^{E}$, where $E$ is the embedding dimension (e.g., here $E=16$). These embeddings are learned through an ``Embedding'' layer and provide a compact, context-rich representation for each district, influencing both the initial latent state and the subsequent latent dynamics.

\subsection{Encoder (GRU)}
The Encoder is responsible for mapping the observed time series data for a given district into an initial latent state $\mathbf{z}_0$. It employs a Gated Recurrent Unit (GRU) to process the sequence \cite{chung2014empirical}. For each time step $t_j$, the input to the GRU is the concatenation of the observed indicator values $\mathbf{x}(t_j)$ and the corresponding district embedding $\mathbf{e}_d$: 
\begin{equation}
    \text{GRU Input}(t_j) = [\mathbf{x}(t_j); \mathbf{e}_d].
\end{equation}

The GRU processes the entire sequence of these concatenated inputs $[\text{GRU Input}(t_0), \dots, \text{GRU Input}(t_{N_t-1})]$ to produce a final hidden state. This hidden state is then passed through a linear layer to obtain the initial latent state $\mathbf{z}_0 \in \mathbb{R}^{L}$, where $L$ is the latent dimension (e.g., here $L=8$). This ensures $\mathbf{z}_0$ captures both the temporal patterns and the specific characteristics of a district.

In particular, a GRU is a recurrent neural network (RNN) module designed to model sequential data while mitigating the vanishing-gradient issue in standard RNNs \cite{cho2014learning}. It maintains a hidden state $h^{(t)}$ that summarizes past information and uses two gates: the update gate $z^{(t)}$ and the reset gate $r^{(t)}$, to regulate information flow when incorporating the current input $x^{(t)}$. Eq.~\eqref{eq:gru_update} below formalizes these operations as \cite{geyer2019deeptma,dey2017gate}:
\begin{align}\label{eq:gru_update}
	z^{(t)} &= \sigma\!\left(W_z x^{(t)} + U_z h^{(t-1)}\right), \nonumber\\
	r^{(t)} &= \sigma\!\left(W_r x^{(t)} + U_r h^{(t-1)}\right), \nonumber\\
	\tilde{h}^{(t)} &= \tanh\!\left(W_h x^{(t)} + U_h \left(r^{(t)} \odot h^{(t-1)}\right)\right), \nonumber\\
	h^{(t)} &= (1 - z^{(t)}) \odot h^{(t-1)} + z^{(t)} \odot \tilde{h}^{(t)},
\end{align}
where $z^{(t)}$ controls the interpolation between the previous hidden state $h^{(t-1)}$ and the candidate state $\tilde{h}^{(t)}$, while $r^{(t)}$ modulates the contribution of $h^{(t-1)}$ in forming $\tilde{h}^{(t-1)}$ via the element-wise product $\odot$.

\begin{figure}[h]
    \centering
    \includegraphics[width=0.8\textwidth]{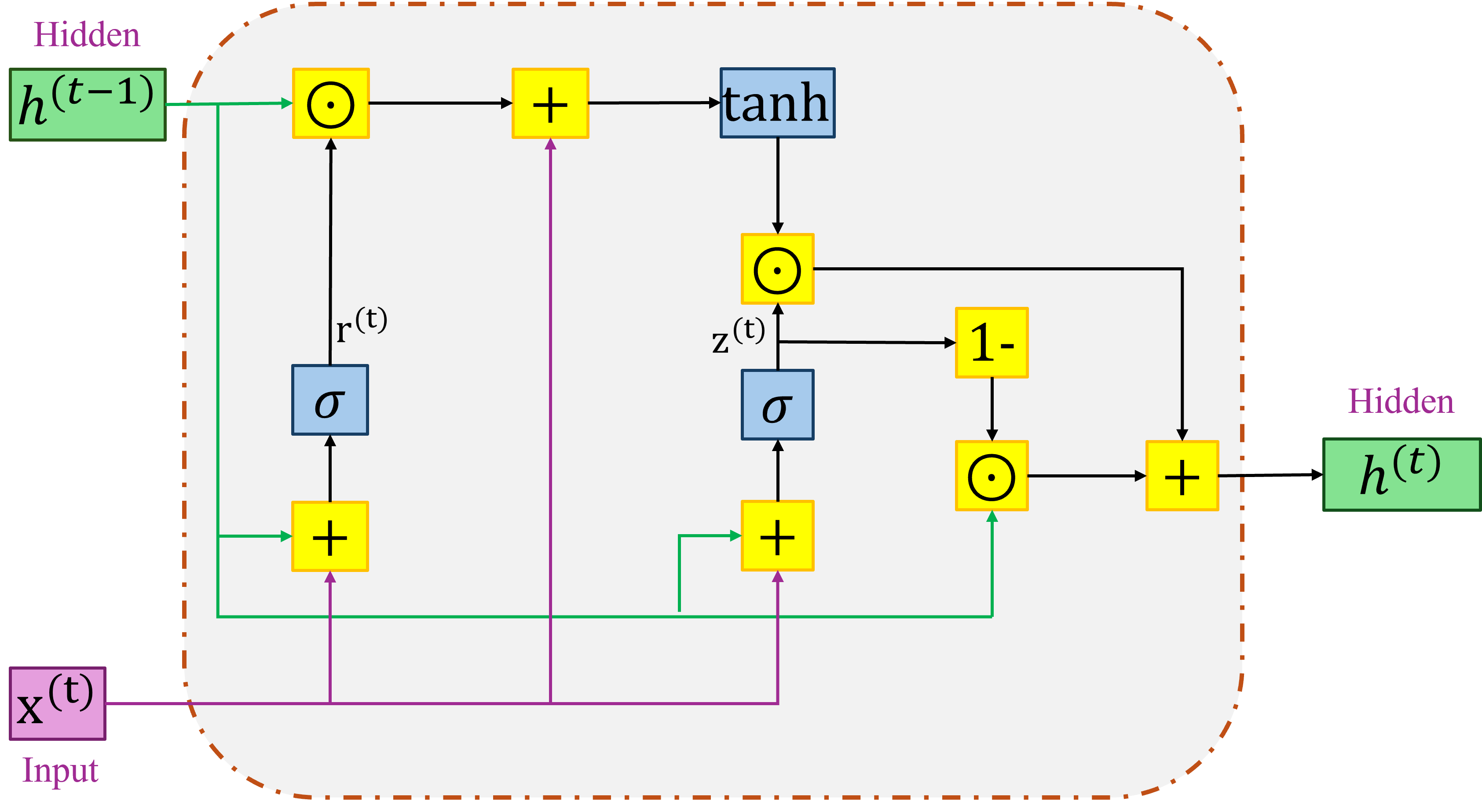}
    \caption{Schematic diagram of a GRU}
    \label{fig:GRU}
\end{figure}

The GRU diagram in Figure~\ref{fig:GRU} provides a graphical interpretation of the same gating mechanism, showing how $x^{(t)}$ and $h^{(t-1)}$ interact through the gates to produce the updated hidden state $h^{(t)}$, which is subsequently mapped to the initial latent representation $\mathbf{z}_0$ for the Neural ODE.

\subsection{District-Conditioned Neural ODE}
This is the core of this Neural ODE model, which defines the continuous latent dynamics. Here, the time derivative of the latent state, $\frac{d\mathbf{z}}{dt}$, is modeled as a neural network $f_{\theta}(\mathbf{z}, \mathbf{e}_d)$. Critically, the district embedding $\mathbf{e}_d$ is concatenated with the current latent state $\mathbf{z}$ at every integration step, making the dynamics district-specific:
\begin{equation}
    \frac{d\mathbf{z}}{dt} = f_{\theta}([\mathbf{z}; \mathbf{e}_d]),
\end{equation}
where function $f_{\theta}$ is implemented as a multi-layer perceptron (MLP) with activation functions (e.g., Tanh) between layers. 

\subsection{ODE Solver}
To obtain the latent trajectory $\mathbf{z}(t)$ at arbitrary query times, the continuous dynamics defined by $f_{\theta}$ are numerically integrated. In this regard, the numerical black box of the 'dopri5' method (Dormand-Prince method of order 5) is used. It is a robust and adaptive-step-size solver. Relative tolerance (rtol) and absolute tolerance (atol) are set to $10^{-3}$ and $10^{-4}$ respectively to control integration accuracy.

\subsection{District-Conditioned Decoder}
The decoder reconstructs the observed indicator values $\hat{\mathbf{x}}(t)$ from the latent state $\mathbf{z}(t)$ and the district embedding $\mathbf{e}_d$. This reconstruction is also district-aware, using the concatenated latent state and embedding as input:
\begin{equation}
    \hat{\mathbf{x}}(t) = g_{\phi}([\mathbf{z}(t); \mathbf{e}_d]),
\end{equation}
where function $g_{\phi}$ is an MLP with a final sigmoid activation function to ensure the output values, representing proportions, lie within the $[0, 1]$ range.

\subsection{Loss Function}
The model is trained by minimizing the Mean Squared Error (MSE) loss between the predicted indicator values $\mathbf{x}_\text{pred}$ and the actual observed values $\mathbf{X}$. The loss is defined as:
\begin{equation}
    \mathcal{L}(\Theta) = \frac{1}{N_d N_t N_i} \sum_{d=1}^{N_d} \sum_{j=1}^{N_t} \sum_{k=1}^{N_i} (x_{\text{pred}}(t_j, d, k) - x(t_j, d, k))^2,
\end{equation}
where $\Theta$ collectively represents all trainable parameters of the model, including those of the Encoder (GRU), the ODE function $f_{\theta}$, the Decoder $g_{\phi}$, and the district embedding matrix $\mathbf{E}$.

\subsection{Training Procedure}
After defining the loss function, the model is trained by backpropagating gradients to update all learnable parameters. In this regard, gradients with respect to the decoder parameters $\phi$ are computed using standard backpropagation, whereas gradients with respect to the Neural ODE parameters $\theta$ are obtained via the adjoint sensitivity method. Gradients for the encoder (GRU) parameters and the district-embedding matrix are again computed using standard backpropagation.

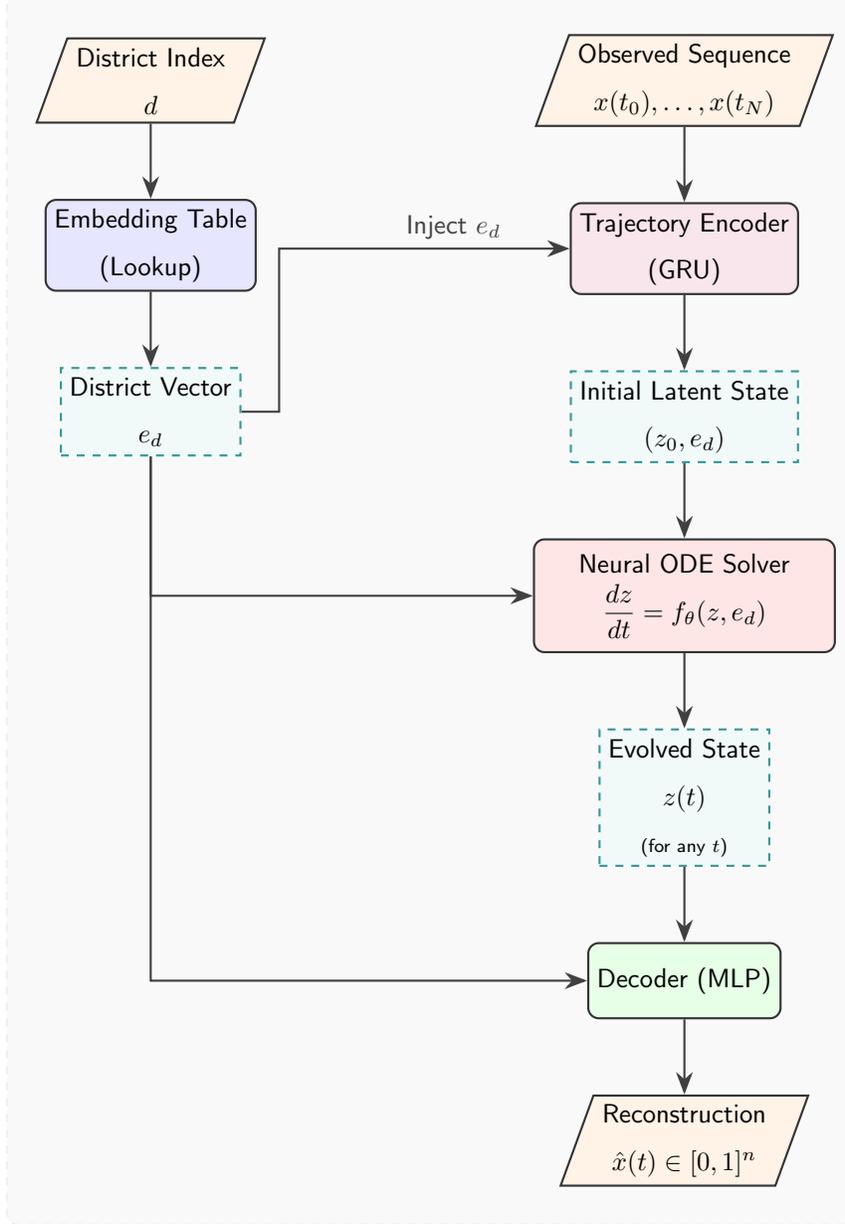
\begin{figure}[h!]
	\centering
	\begin{tikzpicture}[
    node distance=1.5cm and 2cm,
    font=\sffamily,
    % Define styles for different components
    process/.style={rectangle, draw=black!80, fill=blue!10, thick, minimum width=2.5cm, minimum height=1cm, rounded corners, align=center},
    tensor/.style={rectangle, draw=teal!80, fill=teal!5, dashed, thick, minimum width=1.5cm, minimum height=0.8cm, align=center},
    input/.style={trapezium, trapezium left angle=70, trapezium right angle=110, draw=black!80, fill=orange!10, thick, minimum width=1.5cm, align=center},
    arrow/.style={-{Stealth[length=3mm]}, thick, darkgray},
    label/.style={font=\footnotesize\color{gray!90!black}, align=center}
]

    % --- 1. Inputs ---
    \node (idx) [input] {District Index\\$d$};
    \node (seq) [input, right=4cm of idx] {Observed Sequence\\$x(t_0), \dots, x(t_N)$};

    % --- 2. Embedding Layer ---
    \node (embed_table) [process, below=1cm of idx] {Embedding Table\\(Lookup)};
    \node (ed) [tensor, below=1cm of embed_table] {District Vector\\$e_d$};

    % --- 3. Encoder ---
    \node (encoder) [process, below=1cm of seq, fill=purple!10] {Trajectory Encoder\\(GRU)};
    \node (z0) [tensor, below=1.0cm of encoder] {Initial Latent State\\$(z_0,e_d)$};

    % --- 4. Neural ODE Core ---
    % Represents the continuous time evolution
    \node (node_solver) [process, below=1cm of z0, minimum width=4cm, minimum height=1.5cm, fill=red!10] {Neural ODE Solver\\$\displaystyle \frac{dz}{dt} = f_{\theta}(z, e_d)$};
    
    \node (zt) [tensor, below=1.0cm of node_solver] {Evolved State\\$z(t)$ \\ \scriptsize(for any $t$)};

    % --- 5. Decoder ---
    \node (decoder) [process, below=1.0cm of zt, fill=green!10] {Decoder (MLP)};
    \node (output) [input, below=1cm of decoder, shape border rotate=180] {Reconstruction\\$\hat{x}(t) \in [0,1]^{n}$};

    % --- Connections ---

    % District Index -> Embedding -> e_d
    \draw [arrow] (idx) -- (embed_table);
    \draw [arrow] (embed_table) -- (ed);

    % Sequence -> Encoder
    \draw [arrow] (seq) -- (encoder);

    % Injection of e_d into Encoder
    \draw [arrow] (ed.east) -- ++(0.5,0) |- (encoder.west) node[midway, above, pos=0.8] {Inject $e_d$};

    % Encoder -> z0
    \draw [arrow] (encoder) -- (z0);

    % z0 -> Neural ODE
    \draw [arrow] (z0) -- (node_solver);

    % Injection of e_d into Neural ODE
    \draw [arrow] (ed.south) |- ($(node_solver.west) + (-0.5, 0)$) -- (node_solver.west);

    % Neural ODE -> z(t)
    \draw [arrow] (node_solver) -- (zt);

    % z(t) -> Decoder
    \draw [arrow] (zt) -- (decoder);

    % Injection of e_d into Decoder
    \draw [arrow] (ed.south) |- ($(decoder.west) + (-0.5, 0)$) -- (decoder.west);

    % Decoder -> Output
    \draw [arrow] (decoder) -- (output);

    % --- Grouping / Background ---
    \begin{scope}[on background layer]
        \node [fit=(embed_table) (encoder) (output) (ed) (idx) (seq), fill=gray!5, rounded corners, draw=gray!20, dashed, inner sep=0.5cm] (container) {};
    \end{scope}
	\end{tikzpicture}
	\caption{Flow diagram of Neural ODE with district embeddings}
	\label{fig:Flow_diagram_of_Neural_ODE}
\end{figure}

The model is trained using the Adam optimizer with a learning rate (LR) of $10^{-3}$ and a weight decay of $10^{-5}$ to prevent overfitting \cite{kingma2014adam}. A cosine annealing learning rate scheduler is used to adjust the learning rate during training, gradually reducing it over a specified number of epochs ($N_{\text{epochs}}=1000$). This iterative process allows the model to learn the district embeddings, the latent dynamics, and the reconstruction mapping, enabling accurate predictions and forecasts. For visualizing the whole process, Figure~\ref{fig:Flow_diagram_of_Neural_ODE} represents the flow diagram of Neural ODE with district embeddings.

\section{Results and Discussions}\label{sec:results_and_discussions}

This section presents and discusses the results obtained from training and evaluating the District-Embedded Neural ODE model on the poverty indicator time series data. The performance, the insights gained from learned representations, and the forecasting capabilities are analyzed of the model.

\subsection{Model Parameters Setting}

To train the above-described Neural ODE setup for socio-economic dynamics, the GRU encoder, mapping $[\mathbf{x}(t); \mathbf{e}_d] \to \mathbf{z}_0$, is implemented with a single hidden layer of 64 neurons. The latent dynamics function $f_{\theta}:[\mathbf{z}(t), \mathbf{e}_d] \to \mathbf{z}(t)$ is modeled as a neural network with two hidden layers of 64 neurons each, using $Tanh$ activation. The decoder $g_{\phi}:[\mathbf{z}(t); \mathbf{e}_d] \to \mathbf{x}(t)$ is also a neural network with two hidden layers of 64 neurons employing ReLU activation and a Sigmoid function at the output layer.

The district embedding dimension is set to 16, while the latent state dimension is 8. Training is carried out using naive backpropagation for the encoder and decoder, and the adjoint sensitivity method for the Neural ODE component. The model parameters are optimized using the Adam optimizer with a learning rate of 0.001.

\subsection{Loss Convergence}

\begin{figure}[h!]
    \centering
    \includegraphics[width=0.7\textwidth]{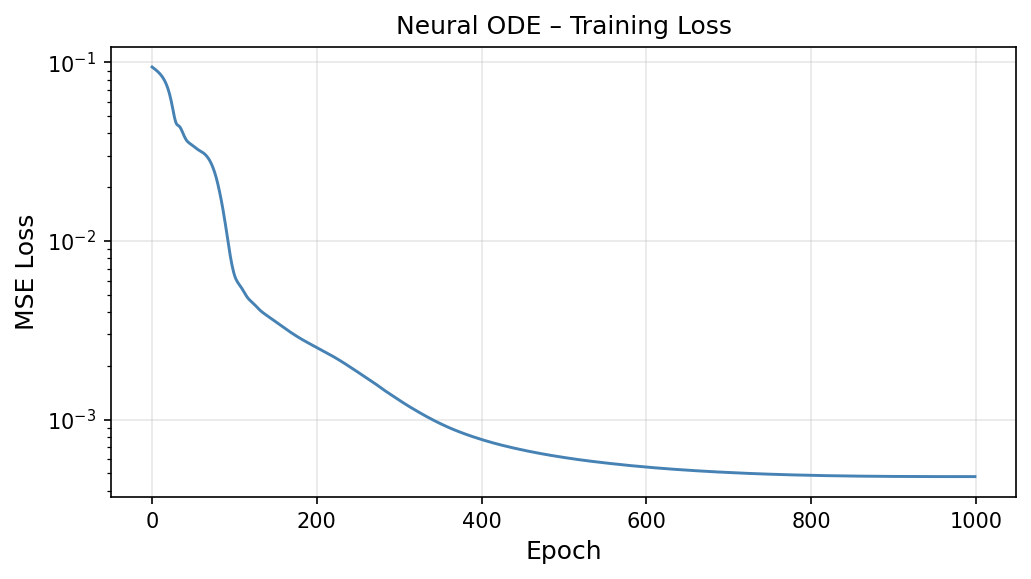}
    \caption{Training Loss Curve: Log-scaled MSE loss over 1000 epochs, demonstrating stable convergence.}
    \label{fig:loss_curve}
\end{figure}

The training process, spanning 1000 epochs, demonstrated stable convergence as evidenced by the training loss curve as depicted in Figure~\ref{fig:loss_curve}. The Mean Squared Error (MSE) loss progressively decreased, reaching a final value of approximately 0.000479, corresponding to a Root Mean Squared Error (RMSE) of 0.021885. This low RMSE indicates that the model effectively learned to reconstruct the observed historical trends of the poverty indicators across all districts.

\subsection{Overall Model Performance}

The Neural ODE model achieved an overall RMSE of 0.021885 on the observed data, demonstrating strong fidelity to the historical trajectories. A detailed breakdown of the RMSE for each indicator highlights the consistent performance of Neural ODE model across different socioeconomic metrics:
\begin{itemize}
    \item Have Access to toilet facility: 0.02203
    \item Use piped drinking water: 0.02818
    \item Use LPG for cooking: 0.02038
    \item Live in a pucca house: 0.02150
    \item Have Electricity connection: 0.02007
    \item Education Secondary or higher (15-49): 0.01774
\end{itemize}
These low individual RMSE values suggest that the model successfully captures the distinct temporal dynamics of each indicator.

\subsection{Learned District Embeddings}

\begin{figure}[h!]
    \centering
    \includegraphics[width=0.8\textwidth]{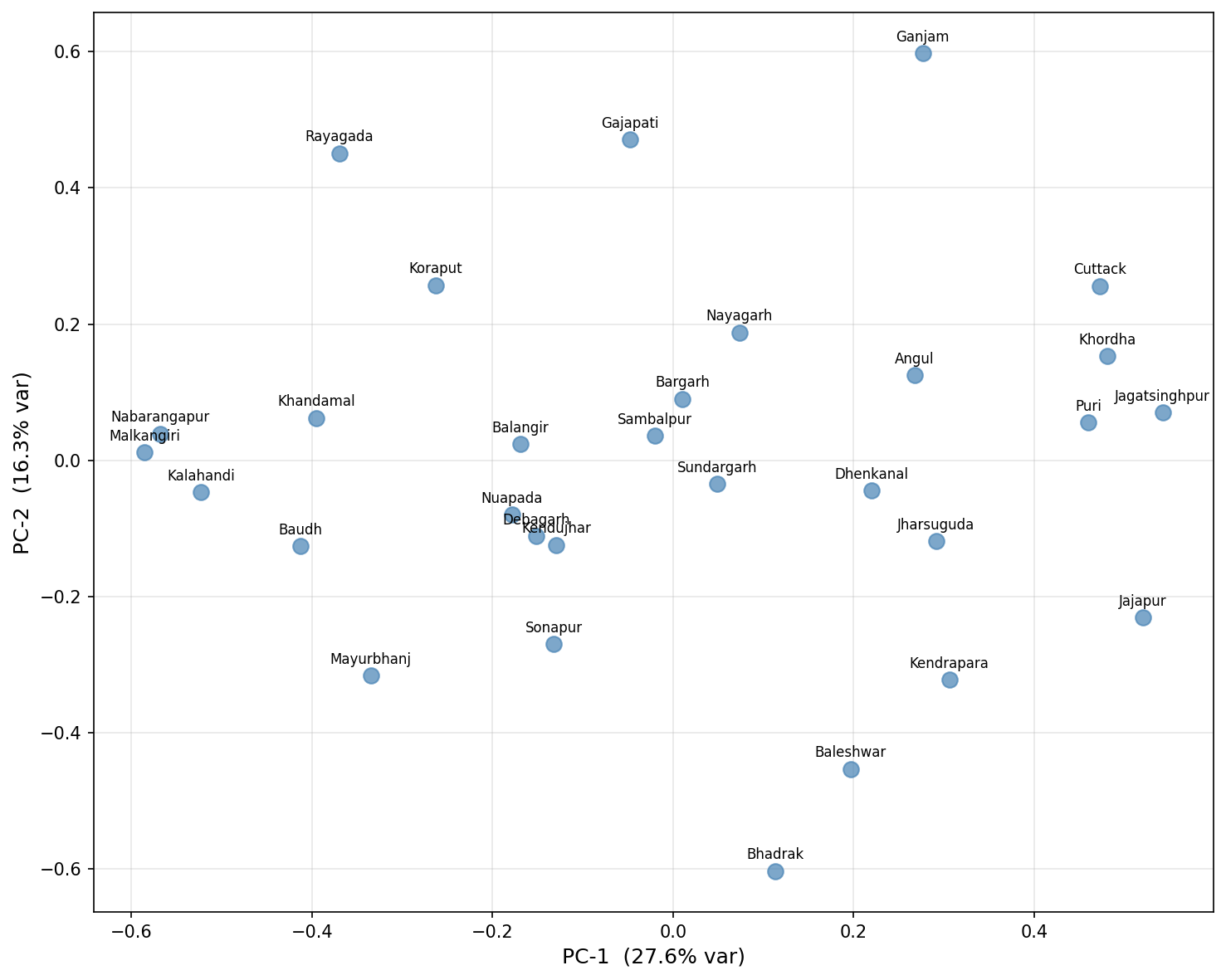}
    \caption{PCA Projection of Learned District Embeddings. Districts with similar poverty dynamics tend to cluster together, reflecting underlying structural groups.}
    \label{fig:embeddings_pca}
\end{figure}

One of the key advantages of the District-Embedded Neural ODE is the learning of unique, low-dimensional embeddings for each district. As illustrated by the PCA projection of these embeddings in Figure~\ref{fig:embeddings_pca}, districts with similar socioeconomic characteristics and trajectories tend to cluster together in the embedding space. This visualization confirms that the model is learning meaningful, context-rich representations for each district, which implicitly encode their unique developmental patterns and influences.

\subsection{Indicator Trajectories and Future Forecasts}

\begin{figure}[h!]
    \centering
    \includegraphics[width=\textwidth]{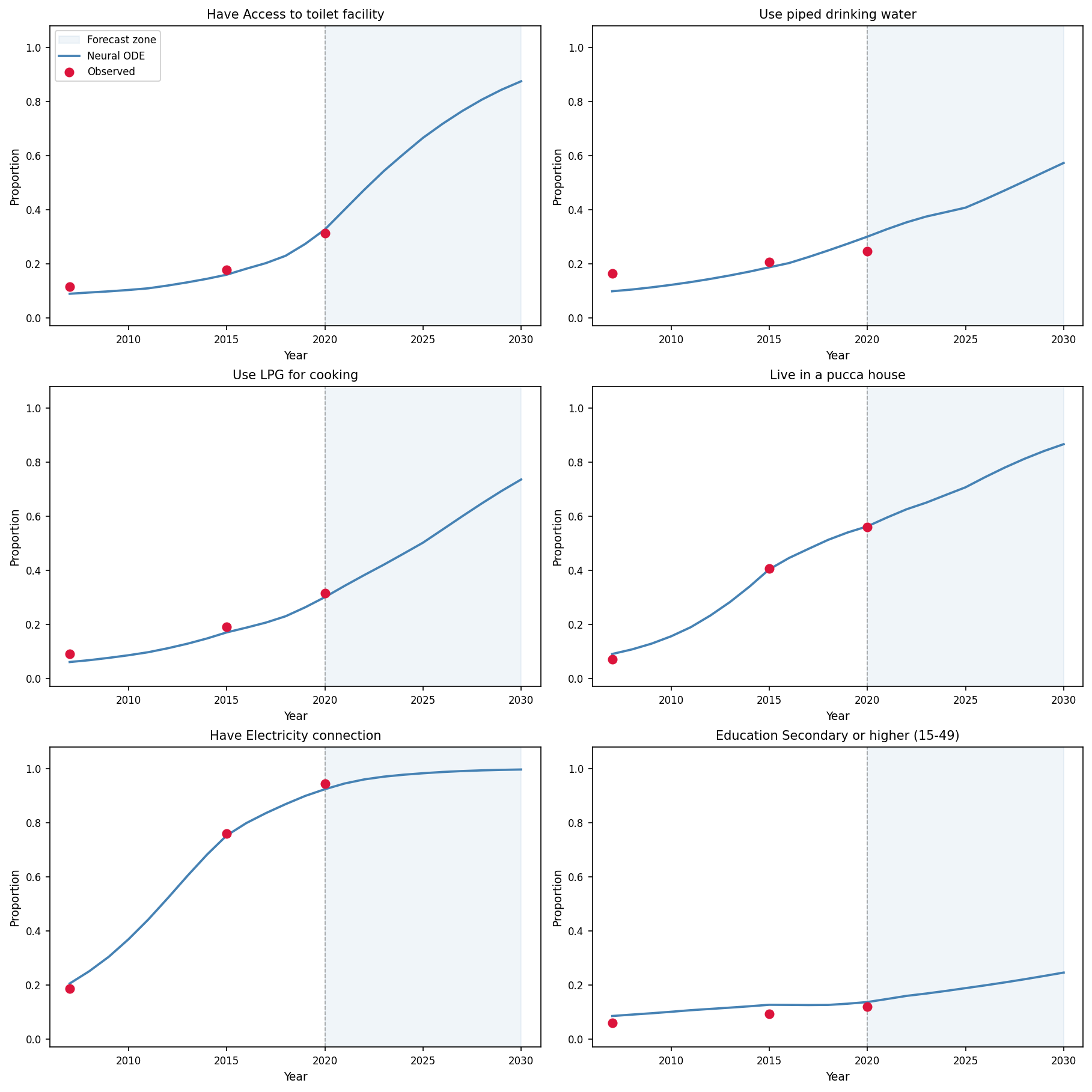}
    \caption{Neural ODE Trajectories for Koraput District; Observed values (red dots), model fit and forecast (steelblue line).}
    \label{fig:koraput_traj}
\end{figure}

\begin{figure}[h!]
    \centering
    \includegraphics[width=\textwidth]{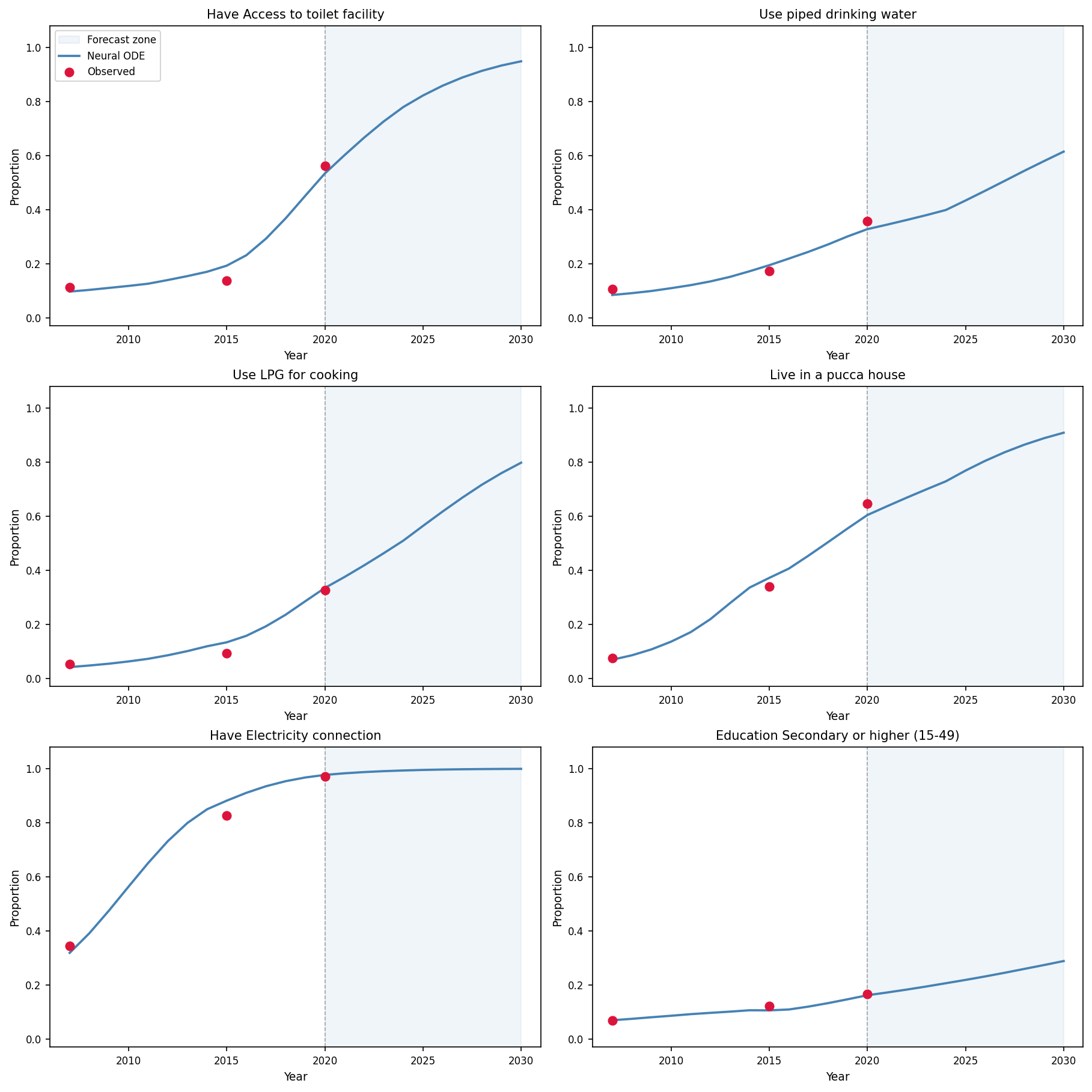}
    \caption{Neural ODE Trajectories for Balangir District; Observed values (red dots), model fit and forecast (steelblue line).}
    \label{fig:balangir_traj}
\end{figure}

\begin{figure}[h!]
    \centering
    \includegraphics[width=\textwidth]{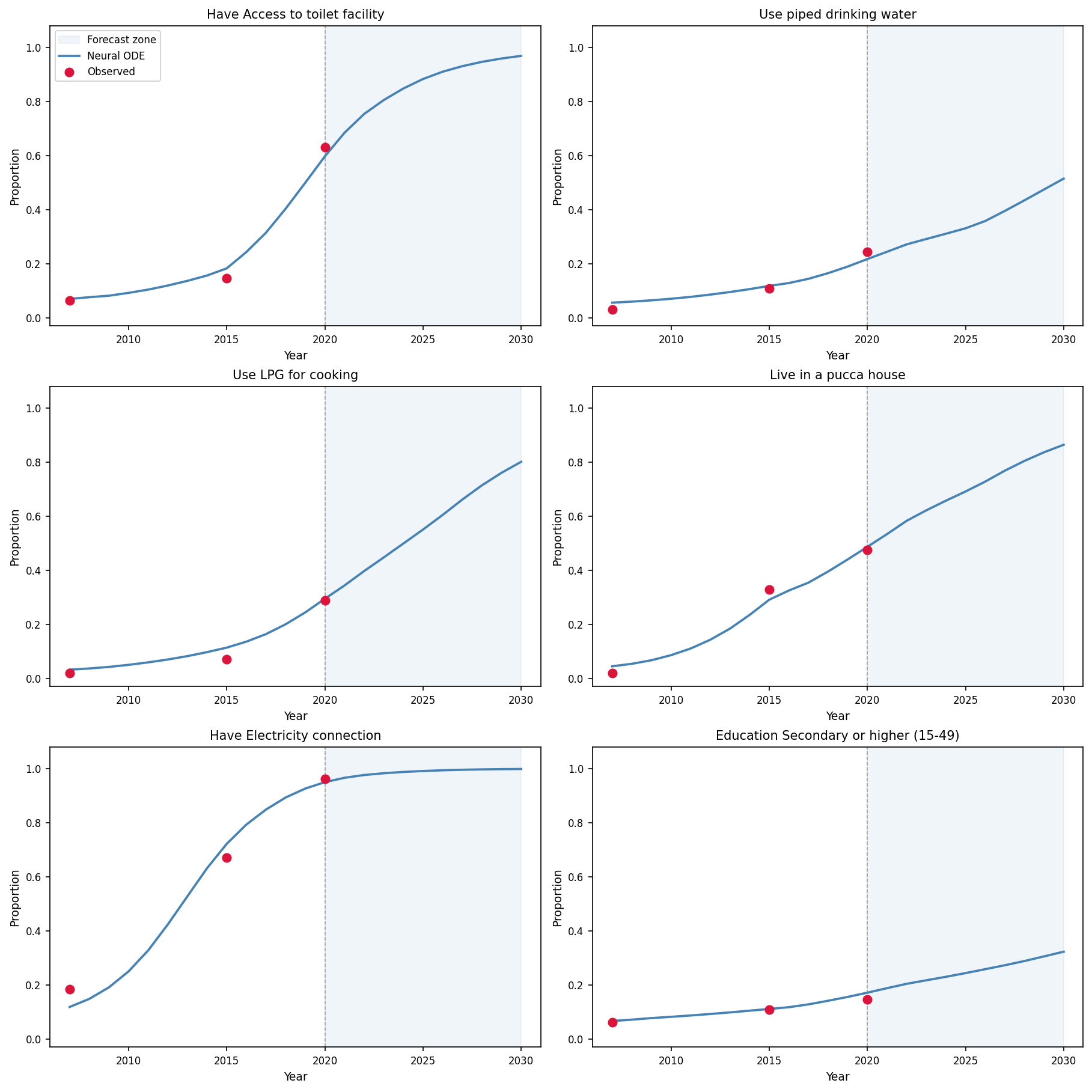}
    \caption{Neural ODE Trajectories for Kalahandi District; Observed values (red dots), model fit and forecast (steelblue line).}
    \label{fig:kalahandi_traj}
\end{figure}

The ability of the Neural ODE to reconstruct observed data and forecast future trends is critical. For visualization of the predicted trajectories, Figures~\ref{fig:koraput_traj}, \ref{fig:balangir_traj}, and \ref{fig:kalahandi_traj} for the KBK (Koraput--Balangir--Kalahandi) districts are presented. These figures show the observed data points (red dots), the fitted trajectories of the Neural ODE over the observed period, and the forecasted trajectories for future years (e.g., 2025 and 2030).

For Koraput, the model accurately captures the historical trends for various indicators. The forecasted trajectories show a continuation of these trends, providing quantitative estimates for future indicator values. For instance, indicators such as 'Have Electricity connection' show a projected increase, while 'Use piped drinking water' also indicates a positive trend. The smooth, continuous nature of the predicted curves highlights the ODE's capacity to model the underlying developmental processes.

In Balangir, the model similarly fits the observed data well. The forecasts suggest specific future trajectories for indicators. Analyzing these plots allows for an understanding of how each indicator is projected to evolve, reflecting the unique context of Balangir as learned by its district embedding. Changes in indicators like 'Education Secondary or higher' are particularly informative for future policy planning.

Kalahandi district indicator trajectories also demonstrate the robustness of the Neural ODE model. The model provides continuous predictions, effectively interpolating between observed time points and extrapolating into the future. These forecasts offer valuable insights into potential future states for each indicator, aiding in proactive strategic development. The consistency of fit of model across diverse districts like Koraput, Balangir, and Kalahandi underscores its generalization capabilities.

\subsection{Forecast Summary}

The model provides quantitative forecasts for future years (e.g., 2025 to 2030) for all districts and indicators. These are presented in a tabular format, offering precise numerical predictions for policymakers and researchers. The ability to predict these values for arbitrary future time points is a direct benefit of the continuous dynamics learned by the Neural ODE.

Tables~\ref{tab:forecast_2026} and~\ref{tab:forecast_2030} summarize the district-wise forecasts for six key poverty-related indicators. The 2026 projections reflect near-term expected levels for access to toilets, piped drinking water, LPG use, pucca housing, electricity, and secondary-or-higher education. These values can be interpreted as predicted proportions, enabling direct comparison across districts and across indicators within a district.

Comparing the 2026 and 2030 forecasts highlights the temporal evolution implied by the learned continuous-time dynamics. Most districts show improvements over time in core infrastructure indicators (toilets, electricity, pucca housing, and LPG), with electricity saturation approaching unity in many cases. The education indicator shows comparatively slower gains, suggesting that improvements in human-capital outcomes may lag behind infrastructure expansion, such differences across indicators provide a useful basis for prioritizing district-specific interventions and monitoring progress over future planning horizons.

\begin{table}[h!]
\centering
\caption{Forecast for Year 2026 across districts and indicators.}
\label{tab:forecast_2026}
\small
\setlength{\tabcolsep}{3pt}
\begin{tabular}{lcccccc}
\hline
\textbf{District} & \textbf{Toilet} & \textbf{Piped water} & \textbf{LPG} & \textbf{Pucca house} & \textbf{Electricity} & \textbf{Education} \\
\hline
Angul & 0.865 & 0.382 & 0.556 & 0.846 & 0.994 & 0.297 \\
Balangir & 0.858 & 0.470 & 0.617 & 0.805 & 0.997 & 0.231 \\
Baleshwar & 0.851 & 0.401 & 0.500 & 0.652 & 0.993 & 0.293 \\
Bargarh & 0.813 & 0.467 & 0.579 & 0.799 & 0.995 & 0.251 \\
Baudh & 0.891 & 0.329 & 0.606 & 0.755 & 0.997 & 0.226 \\
Bhadrak & 0.852 & 0.248 & 0.475 & 0.640 & 0.995 & 0.306 \\
Cuttack & 0.853 & 0.608 & 0.727 & 0.911 & 0.995 & 0.340 \\
Debagarh & 0.844 & 0.329 & 0.460 & 0.661 & 0.992 & 0.259 \\
Dhenkanal & 0.786 & 0.346 & 0.571 & 0.798 & 0.993 & 0.283 \\
Gajapati & 0.810 & 0.567 & 0.601 & 0.851 & 0.994 & 0.198 \\
Ganjam & 0.890 & 0.771 & 0.844 & 0.956 & 0.998 & 0.293 \\
Jagatsinghpur & 0.863 & 0.434 & 0.597 & 0.891 & 0.997 & 0.335 \\
Jajapur & 0.807 & 0.276 & 0.525 & 0.821 & 0.994 & 0.330 \\
Jharsuguda & 0.845 & 0.475 & 0.672 & 0.830 & 0.992 & 0.343 \\
Kalahandi & 0.910 & 0.358 & 0.604 & 0.728 & 0.994 & 0.258 \\
Kendrapara & 0.844 & 0.357 & 0.547 & 0.787 & 0.996 & 0.318 \\
Kendujhar & 0.737 & 0.391 & 0.504 & 0.652 & 0.982 & 0.264 \\
Khandamal & 0.901 & 0.304 & 0.556 & 0.788 & 0.996 & 0.247 \\
Khordha & 0.909 & 0.702 & 0.850 & 0.940 & 0.997 & 0.414 \\
Koraput & 0.717 & 0.438 & 0.551 & 0.744 & 0.988 & 0.198 \\
Malkangiri & 0.860 & 0.347 & 0.490 & 0.684 & 0.997 & 0.172 \\
Mayurbhanj & 0.771 & 0.335 & 0.384 & 0.484 & 0.982 & 0.245 \\
Nabarangapur & 0.827 & 0.365 & 0.460 & 0.644 & 0.992 & 0.186 \\
Nayagarh & 0.908 & 0.536 & 0.749 & 0.911 & 0.999 & 0.270 \\
Nuapada & 0.830 & 0.310 & 0.505 & 0.741 & 0.995 & 0.211 \\
Puri & 0.916 & 0.506 & 0.696 & 0.922 & 0.998 & 0.348 \\
Rayagada & 0.856 & 0.610 & 0.626 & 0.813 & 0.995 & 0.192 \\
Sambalpur & 0.836 & 0.573 & 0.660 & 0.778 & 0.991 & 0.297 \\
Sonapur & 0.885 & 0.373 & 0.682 & 0.829 & 0.998 & 0.277 \\
Sundargarh & 0.856 & 0.529 & 0.661 & 0.790 & 0.991 & 0.329 \\
\hline
\end{tabular}
\end{table}

\begin{table}[h!]
\centering
\caption{Forecast for Year 2030 across districts and indicators.}
\label{tab:forecast_2030}
\small
\setlength{\tabcolsep}{3pt}
\begin{tabular}{lcccccc}
\hline
\textbf{District} & \textbf{Toilet} & \textbf{Piped water} & \textbf{LPG} & \textbf{Pucca house} & \textbf{Electricity} & \textbf{Education} \\
\hline
Angul & 0.945 & 0.524 & 0.748 & 0.923 & 0.998 & 0.362 \\
Balangir & 0.949 & 0.615 & 0.798 & 0.909 & 0.999 & 0.288 \\
Baleshwar & 0.931 & 0.521 & 0.680 & 0.792 & 0.998 & 0.346 \\
Bargarh & 0.925 & 0.603 & 0.762 & 0.901 & 0.999 & 0.306 \\
Baudh & 0.962 & 0.483 & 0.800 & 0.884 & 0.999 & 0.289 \\
Bhadrak & 0.938 & 0.360 & 0.671 & 0.792 & 0.999 & 0.365 \\
Cuttack & 0.936 & 0.718 & 0.853 & 0.954 & 0.999 & 0.398 \\
Debagarh & 0.938 & 0.464 & 0.670 & 0.811 & 0.998 & 0.317 \\
Dhenkanal & 0.905 & 0.479 & 0.751 & 0.892 & 0.998 & 0.343 \\
Gajapati & 0.924 & 0.702 & 0.785 & 0.929 & 0.999 & 0.248 \\
Ganjam & 0.958 & 0.848 & 0.925 & 0.980 & 1.000 & 0.348 \\
Jagatsinghpur & 0.944 & 0.577 & 0.780 & 0.946 & 0.999 & 0.405 \\
Jajapur & 0.914 & 0.398 & 0.713 & 0.906 & 0.998 & 0.395 \\
Jharsuguda & 0.931 & 0.602 & 0.814 & 0.910 & 0.998 & 0.405 \\
Kalahandi & 0.969 & 0.515 & 0.801 & 0.864 & 0.999 & 0.323 \\
Kendrapara & 0.935 & 0.493 & 0.736 & 0.889 & 0.999 & 0.381 \\
Kendujhar & 0.874 & 0.512 & 0.681 & 0.794 & 0.995 & 0.315 \\
Khandamal & 0.967 & 0.460 & 0.776 & 0.900 & 0.999 & 0.310 \\
Khordha & 0.962 & 0.787 & 0.922 & 0.970 & 0.999 & 0.469 \\
Koraput & 0.875 & 0.573 & 0.736 & 0.866 & 0.997 & 0.245 \\
Malkangiri & 0.953 & 0.488 & 0.717 & 0.834 & 0.999 & 0.222 \\
Mayurbhanj & 0.903 & 0.413 & 0.557 & 0.621 & 0.994 & 0.296 \\
Nabarangapur & 0.939 & 0.477 & 0.670 & 0.787 & 0.998 & 0.237 \\
Nayagarh & 0.968 & 0.684 & 0.884 & 0.961 & 1.000 & 0.337 \\
Nuapada & 0.935 & 0.450 & 0.719 & 0.870 & 0.999 & 0.265 \\
Puri & 0.970 & 0.655 & 0.851 & 0.965 & 1.000 & 0.422 \\
Rayagada & 0.947 & 0.735 & 0.804 & 0.911 & 0.999 & 0.239 \\
Sambalpur & 0.927 & 0.684 & 0.805 & 0.882 & 0.998 & 0.350 \\
Sonapur & 0.959 & 0.525 & 0.842 & 0.921 & 1.000 & 0.345 \\
Sundargarh & 0.936 & 0.649 & 0.807 & 0.889 & 0.997 & 0.387 \\
\hline
\end{tabular}
\end{table}

\section{Conclusion}\label{sec:conclusion}

This chapter demonstrates the utility of a district-embedded Neural ODE framework for modeling poverty-related indicators as a continuous-time dynamical system. By combining a GRU-based encoder with district embeddings and a neural ODE dynamics function, the proposed approach captures nonlinear temporal patterns, accommodates irregular observation times, and provides smooth trajectories for both interpolation and forecasting. The empirical results indicate that the learned representations help differentiate district-specific dynamics, while the fitted and forecasted trajectories offer an interpretable summary of expected indicator evolution over time.

There are some limitations also of this study which are explained below. Firsty, the analysis relies on a relatively short time span (2007--2020) which may introduce bias and can lead to overconfident extrapolations, especially for longer horizons. Secondly, the model is purely data-driven and does not explicitly incorporate exogenous drivers (e.g., policy interventions, migration, shocks, or macroeconomic variables), which can substantially affect poverty trajectories and may reduce forecast reliability under structural change. Thirdly, uncertainty quantification is limited, as point forecasts are reported without calibrated predictive intervals. Future work could incorporate Bayesian or ensemble approaches to better characterize uncertainty and improve robustness for policy use.

%\section*{Statements and Declarations}
%No funding was received to assist with the preparation of this manuscript.
%\begin{itemize} 
%\item \textbf{Data:} Data sharing not applicable to this article as no datasets were generated or analysed during the current study.
%\item \textbf{Funding:} No funding was received to assist with the preparation of this manuscript.
%\item \textbf{Competing interests:} The authors have no competing interests to declare that are relevant to the content of this article.
%\end{itemize} 

\section*{Acknowledgement}
The first author would like to acknowledge the Poverty Alleviation Research Centre (PARC) at the National Institute of Technology Rourkela for providing fellowship support to carry out this research.

\bibliography{reff}

@inproceedings{he2016deep,
  title={Deep Residual Learning for Image Recognition},
  author={He, Kaiming and Zhang, Xiangyu and Ren, Shaoqing and Sun, Jian},
  booktitle={Proceedings of the IEEE Conference on Computer Vision and Pattern Recognition (CVPR)},
  pages={770--778},
  year={2016}
}

@book{pontryagin2018mathematical,
  title={The Mathematical Theory of Optimal Processes},
  author={Pontryagin, Lev S. and Boltyanskii, V. G. and Gamkrelidze, R. V. and Mishchenko, E. F.},
  year={2018},
  publisher={Routledge}
}

@article{carr2009evaluating,
  title={Evaluating poverty--environment dynamics},
  author={Carr, Edward R and Kettle, NP and Hoskins, A},
  journal={International Journal of Sustainable Development \& World Ecology},
  volume={16},
  number={2},
  pages={87--93},
  year={2009},
  publisher={Taylor \& Francis}
}

@article{chen2018neural,
  title={Neural ordinary differential equations},
  author={Chen, Ricky TQ and Rubanova, Yulia and Bettencourt, Jesse and Duvenaud, David K},
  journal={Advances in neural information processing systems},
  volume={31},
  year={2018}
}

@article{rubanova2019latent,
  title={Latent ordinary differential equations for irregularly-sampled time series},
  author={Rubanova, Yulia and Chen, Ricky TQ and Duvenaud, David K},
  journal={Advances in neural information processing systems},
  volume={32},
  year={2019}
}

@article{grathwohl2018ffjord,
  title={Ffjord: Free-form continuous dynamics for scalable reversible generative models},
  author={Grathwohl, Will and Chen, Ricky TQ and Bettencourt, Jesse and Sutskever, Ilya and Duvenaud, David},
  journal={arXiv preprint arXiv:1810.01367},
  year={2018}
}

@article{oh2025comprehensive,
  title={Comprehensive review of neural differential equations for time series analysis},
  author={Oh, YongKyung and Kam, Seungsu and Lee, Jonghun and Lim, Dong-Young and Kim, Sungil and Bui, Alex},
  journal={arXiv preprint arXiv:2502.09885},
  year={2025}
}

@inproceedings{geyer2019deeptma,
  title={DeepTMA: Predicting effective contention models for network calculus using graph neural networks},
  author={Geyer, Fabien and Bondorf, Steffen},
  booktitle={IEEE INFOCOM 2019-IEEE Conference on Computer Communications},
  pages={1009--1017},
  year={2019},
  organization={IEEE}
}

@inproceedings{cho2014learning,
  title={Learning phrase representations using RNN encoder--decoder for statistical machine translation},
  author={Cho, Kyunghyun and Van Merri{\"e}nboer, Bart and Gul{\c{c}}ehre, {\c{C}}a{\u{g}}lar and Bahdanau, Dzmitry and Bougares, Fethi and Schwenk, Holger and Bengio, Yoshua},
  booktitle={Proceedings of the 2014 conference on empirical methods in natural language processing (EMNLP)},
  pages={1724--1734},
  year={2014}
}

@article{chung2014empirical,
  title={Empirical evaluation of gated recurrent neural networks on sequence modeling},
  author={Chung, Junyoung and Gulcehre, Caglar and Cho, KyungHyun and Bengio, Yoshua},
  journal={arXiv preprint arXiv:1412.3555},
  year={2014}
}

@article{chien2022learning,
  title={Learning continuous-time dynamics with attention},
  author={Chien, Jen-Tzung and Chen, Yi-Hsiang},
  journal={IEEE Transactions on Pattern Analysis and Machine Intelligence},
  volume={45},
  number={2},
  pages={1906--1918},
  year={2022},
  publisher={IEEE}
}

@article{kumar2025assessing,
  title={Assessing changes in wealth index using primary survey data},
  author={Kumar, Sandeep and Chakraverty, S and Sethi, Narayan},
  journal={Socio-Economic Planning Sciences},
  volume={98},
  pages={102115},
  year={2025},
  publisher={Pergamon}
}

@article{kumar2023multidimensional,
  title={Multidimensional poverty: CMPI development, spatial analysis and clustering},
  author={Kumar, Sandeep and Chakraverty, Snehashish and Sethi, Narayan},
  journal={Social Indicators Research},
  volume={169},
  number={1},
  pages={647--670},
  year={2023},
  publisher={Springer Netherlands Dordrecht}
}

@article{deaton1985panel,
  title={Panel data from time series of cross-sections},
  author={Deaton, Angus},
  journal={Journal of econometrics},
  volume={30},
  number={1-2},
  pages={109--126},
  year={1985},
  publisher={Elsevier}
}

@book{pesaran2015time,
  title={Time series and panel data econometrics},
  author={Pesaran, M Hashem},
  year={2015},
  publisher={Oxford University Press}
}

@article{wu2015economic,
  title={Economic development and socioeconomic inequality of well-being: A cross-sectional time-series analysis of urban China, 2003--2011},
  author={Wu, Hania Fei and Tam, Tony},
  journal={Social Indicators Research},
  volume={124},
  number={2},
  pages={401--425},
  year={2015},
  publisher={Springer}
}

@article{steele2017mapping,
  title={Mapping poverty using mobile phone and satellite data},
  author={Steele, Jessica E and Sunds{\o}y, P{\aa}l Roe and Pezzulo, Carla and Alegana, Victor A and Bird, Tomas J and Blumenstock, Joshua and Bjelland, Johannes and Eng{\o}-Monsen, Kenth and De Montjoye, Yves-Alexandre and Iqbal, Asif M and others},
  journal={Journal of The Royal Society Interface},
  volume={14},
  number={127},
  year={2017},
  publisher={The Royal Society}
}

@article{jean2016combining,
  title={Combining satellite imagery and machine learning to predict poverty},
  author={Jean, Neal and Burke, Marshall and Xie, Michael and Alampay Davis, W Matthew and Lobell, David B and Ermon, Stefano},
  journal={Science},
  volume={353},
  number={6301},
  pages={790--794},
  year={2016},
  publisher={American Association for the Advancement of Science}
}

@article{khoun2025mapping,
  title={Mapping the dimensions of poverty through big data, socioeconomic surveys and machine learning in Cambodia},
  author={Khoun, Theara and Poortinga, Ate and Thwal, Nyein Soe and Gonz{\'a}lez de Alba, Iv{\'a}n and McMahon, Andrea and Mendez, Carlos},
  journal={Social Indicators Research},
  volume={180},
  number={3},
  pages={1593--1618},
  year={2025},
  publisher={Springer}
}

@article{mcdermott2019bayesian,
  title={Bayesian recurrent neural network models for forecasting and quantifying uncertainty in spatial-temporal data},
  author={McDermott, Patrick L and Wikle, Christopher K},
  journal={Entropy},
  volume={21},
  number={2},
  pages={184},
  year={2019},
  publisher={MDPI}
}

@article{emshagin2022short,
  title={Short-term prediction of household electricity consumption using customized LSTM and GRU models},
  author={Emshagin, Saad and Halim, Wayes Koroni and Kashef, Rasha},
  journal={arXiv preprint arXiv:2212.08757},
  year={2022}
}

@article{kidger2020neural,
  title={Neural controlled differential equations for irregular time series},
  author={Kidger, Patrick and Morrill, James and Foster, James and Lyons, Terry},
  journal={Advances in neural information processing systems},
  volume={33},
  pages={6696--6707},
  year={2020}
}

@article{kidger2022neural,
  title={On neural differential equations},
  author={Kidger, Patrick},
  journal={arXiv preprint arXiv:2202.02435},
  year={2022}
}

@article{paszke2017automatic,
  title={Automatic differentiation in pytorch},
  author={Paszke, Adam and Gross, Sam and Chintala, Soumith and Chanan, Gregory and Yang, Edward and DeVito, Zachary and Lin, Zeming and Desmaison, Alban and Antiga, Luca and Lerer, Adam},
  year={2017}
}

@online{dhs_dataset,
	author  = {{The DHS Program}},
	title   = {NFHS-4 and NFHS-5 survey raw data by DHS program},
	year    = {2020},
	howpublished = {\url{https://dhsprogram.com/data/available-datasets.cfm}},
	note = {Accessed: 2026-02-25}
}

@online{dlhs3_dataset,
	author  = {{Government of India}},
	title   = {Population and household characteristics, DLHS-III},
	year    = {2014},
	howpublished = {\url{https://www.data.gov.in/resource/population-and-household-characteristics-total-dlhs-iii}},
	note = {Accessed: 2026-02-25}
}

@inproceedings{dey2017gate,
  title={Gate-variants of gated recurrent unit (GRU) neural networks},
  author={Dey, Rahul and Salem, Fathi M},
  booktitle={2017 IEEE 60th international midwest symposium on circuits and systems (MWSCAS)},
  pages={1597--1600},
  year={2017},
  organization={IEEE}
}

@article{kingma2014adam,
  title={Adam: A method for stochastic optimization},
  author={Kingma, Diederik P and Ba, Jimmy},
  journal={arXiv preprint arXiv:1412.6980},
  year={2014}
}

@article{sahu2026physics,
  title={Physics-informed functional link with theory of functional connections technique for solving differential equations},
  author={Sahu, Iswari and Kumar, Sandeep and Chakraverty, S},
  journal={Neurocomputing},
  pages={132795},
  year={2026},
  publisher={Elsevier}
}

@article{kumar2025comparative,
  title={A comparative study of Center--Radius and Lower--Upper type interval neural network methods in uncertainty modeling},
  author={Kumar, Sandeep and Chakraverty, Snehashish and Sethi, Narayan},
  journal={Applied Soft Computing},
  volume={180},
  pages={113347},
  year={2025},
  publisher={Elsevier}
}

@inproceedings{kumar2023physics,
  title={Physics-informed machine learning framework for approximating the modified degasperis-procesi equation},
  author={Kumar, Sandeep and Sahoo, Arup Kumar and Chakraverty, S},
  booktitle={2023 International Conference on Ambient Intelligence, Knowledge Informatics and Industrial Electronics (AIKIIE)},
  pages={1--6},
  year={2023},
  organization={IEEE}
}

@article{raissi2019physics,
  title={Physics-informed neural networks: A deep learning framework for solving forward and inverse problems involving nonlinear partial differential equations},
  author={Raissi, Maziar and Perdikaris, Paris and Karniadakis, George E},
  journal={Journal of Computational physics},
  volume={378},
  pages={686--707},
  year={2019},
  publisher={Elsevier}
}

%\textbf{\textbf{\textbf{}}}
\end{document}